\documentclass[a4paper,fleqn]{cas-sc}

\usepackage[authoryear,longnamesfirst]{natbib}
\usepackage{graphicx}
\usepackage{subfigure}
\usepackage{listings}
\usepackage{xcolor}
\usepackage{amssymb}
\usepackage{amsmath}
\usepackage{multirow}
\usepackage{array}
\usepackage{algorithm}
\usepackage{algpseudocode}
\usepackage{lineno,hyperref}
\usepackage{booktabs}
\usepackage{float}

\def\tsc#1{\csdef{#1}{\textsc{\lowercase{#1}}\xspace}}
\tsc{ABM}
\tsc{EGT}
\tsc{NLP}
\tsc{LLM}

\begin{document}
\let\WriteBookmarks\relax
\def\floatpagepagefraction{1}
\def\textpagefraction{.001}

\shorttitle{Bridging Semantics and Behavior}

\shortauthors{Blekanov I., Gubar E., Fun X.}

\title [mode = title]{Bridging Semantics and Behavior: An Agent-Based Evolutionary Model of Topic Competition with BERTopic}

\author[1]{Blekanov I.}[   role=Researcher,
    orcid=0000-0002-7305-1429]
    
\author[1]{Gubar E.}[
    role=Researcher,
    orcid=0000-0002-8970-1617]


\cormark[2]

\ead{alyona.gubar@gmail.com, e.gubar@spbu.ru}

\affiliation[1]{organization={St. Petersburg State University},
    addressline={7/9 Universitetskaya nab., },
    city={St. Petersburg},
    postcode={ 199034},
    country={Russia}}

\cortext[2]{Gubar E.}
\author[1]{Fan X.}[role=Researcher, ]
\credit{Blekanov V. Conceptualization, Methodology, Software, Validation, Formal analysis, Writing - Original draft, Writing - Review \& Editing Gubar E. -Conceptualization, Methodology,  Validation, Formal analysis, Writing - Original draft, Writing - Review \& Editing, Fan X.   Software, Validation, Formal analysis, Writing - Original draft.}

\begin{abstract}
Understanding the psychological drivers of information diffusion and opinion evolution in social media is a central challenge in computational social science. While users are constantly exposed to competing topics, their decisions to engage with and propagate information are influenced by social learning, perceived utility, and bounded rationality. Traditional models often decouple semantic content from behavioral dynamics, limiting their ability to explain how individual cognition and social interaction jointly shape large-scale discourse. To address this gap, we propose a novel language-independent computational framework that integrates semantic topic modeling with evolutionary game theory and agent-based simulation to model how individuals adopt and propagate topics in social networks. Using BERTopic, we first extract semantically coherent topics from a large-scale Weibo dataset collected using openly available API tools. We then quantify user engagement through a unified influence metric based on likes, comments, and reposts. In the NetLogo environment, agents update their topic preferences through local imitation, guided by payoff comparisons and behavioral noise—reflecting real-world cognitive constraints and imperfect decision-making. Experiments conducted on two real-world datasets show that the model successfully replicates observed topic diffusion trends, highlighting the role of social learning frequency and payoff sensitivity in shaping public attention. Our findings illustrate how micro-level psychological mechanisms—such as imitation and noise-tolerant choice—can lead to emergent macro-level discourse patterns. This study contributes to the interdisciplinary understanding of human-computer interaction by bridging semantic analysis with behavioral simulation, offering insights into the cognitive and social foundations of online collective behavior. The framework also provides a basis for applications such as misinformation resilience, opinion dynamics modeling, and adaptive content curation.
\end{abstract}




\maketitle

\section{Introduction}
The architecture of modern public discourse has been fundamentally reshaped by social media platforms like Weibo, X (formerly Twitter), and Facebook. These platforms have democratized information production and dissemination, transforming them into primary arenas for public opinion formation \cite{pentina2014}. However, this democratization has given rise to a complex and often unpredictable ecosystem of information diffusion. A central puzzle in computational social science is why certain topics achieve viral status, capturing massive public attention, while others, seemingly of equal importance, fade into obscurity \cite{lazer2020}. This selective attention and the resultant dynamics of topic competition represent a critical challenge for modeling and predicting large-scale social phenomena.

Traditional approaches to understanding these dynamics often fall into two distinct, and largely isolated, paradigms. The first, rooted in computational linguistics and natural language processing (NLP), focuses on content. It employs techniques like topic modeling to extract latent semantic structures from vast corpora of user-generated text, aiming to understand "what" is being discussed \cite{blei2012, bodrunova2019, grootendorst2022}. The second paradigm, emerging from network science and agent-based modeling (ABM), focuses on structure and interaction. It leverages data on user relationships and engagements (likes, shares, replies) to model "how" information propagates through the social fabric, simulating behavioral contagion and influence cascades \cite{helbing2012, lerman2010}. While both paradigms offer powerful insights, their decoupling limits a holistic understanding. A comprehensive model must bridge this gap, connecting the semantic meaning of discourse with the behavioral mechanisms that drive its spread and evolution \cite{bodrunova2019, bodrunova2023}.

This gap is not merely methodological but also theoretical. Classical models of information diffusion and opinion formation often assume a high degree of individual rationality, where users make optimal choices based on complete information. However, a wealth of psychological research contradicts this assumption, demonstrating that human decision-making is characterized by bounded rationality \cite{simon1990}. In the context of social media, users are faced with information overload, attention scarcity, and complex social signals. Their decisions—to engage with a post, adopt a topic, or share information—are not perfectly rational optimizations but are instead guided by cognitive shortcuts, such as social learning and imitation \cite{bandura1971, rendell2010}. Users often look to the behaviors and expressed opinions of others in their network as heuristic cues for what is important or valid, a process amplified by platform algorithms that highlight popular content \cite{bakshy2015}. This imitation is further perturbed by behavioral noise, reflecting individual differences, emotional states, and the inherent uncertainty of the online environment~\cite{lorits2023evolutionary}.

This perspective aligns with the emerging critique of the "deliberative" ideal of public discourse, as outlined in the project description. The broader project of this research on the concept of "cumulative deliberation" posits that online consensus is not typically reached through rational, round-table discussion but through the accumulation of individual expressions, often devoid of deliberative intent \cite{bodrunova2023}. This view re-centers the focus on the micro-level mechanisms—imitation, bounded rationality, and noise-tolerant choice—that lead to the emergent, macro-level patterns of topic dominance and decay. It acknowledges the user's right to be "non-deliberative" in their natural online behavior, while recognizing that the aggregate of these behaviors holds significant socio-political weight. This reality is reflected in the dissonant and antagonistic nature of contemporary online spaces, where rational dialogue is often supplanted by emotional expression and the constant instability of actors who rapidly enter and exit discussions \cite{smoliarova2021, carpentier2017}. The very emotional architecture of social media platforms incentivizes affective, rather than purely rational, engagement, further distancing user behavior from classical deliberative norms \cite{wahljorgensen2018}. In this environment, opinions can be seen as forming not just through linear argumentation but through the accumulation and competition of discursive structures and ideologies that frame issues in particular ways \cite{vandijk2000}.

To capture this complexity, we propose a novel, language-independent computational framework that integrates semantic topic modeling with evolutionary game theory (EGT) and agent-based simulation. Our work makes several key contributions:

\begin{enumerate}
\item Theoretical Integration. We bridge the semantic-content and behavioral-interaction paradigms by formally treating topics, extracted via advanced NLP, as strategic options in a social game. This allows us to model how the meaning of information influences its competitive fitness in the attention economy.
\item Psychological Plausibility. Grounded in the principles of bounded rationality and social learning, our agent-based model in NetLogo operationalizes micro-level decision-making. Agents update their topic preferences through local imitation, guided by payoff comparisons (reflecting perceived social utility) and behavioral noise (reflecting cognitive constraints and imperfect judgment). This directly addresses the Computers in Human Behavior scope by focusing on how computer-mediated environments shape human learning and social interaction.
\item Methodological Pipeline. We establish a reproducible pipeline (framework) from raw text to behavioral simulation. Using BERTopic, we extract coherent topics from a large-scale Weibo dataset. We then develop a novel method to transform user engagement metrics (likes, comments, reposts) into a topic payoff matrix, which drives the evolutionary dynamics within the agent-based model.
\item Empirical Validation and Insight. Through rigorous experiments on two real-world datasets centered on climate change discussions, we demonstrate our model's ability to replicate observed macro-level diffusion trends. We quantitatively explore the impact of key behavioral parameters—such as social learning frequency (prob-revision) and payoff sensitivity ($\alpha$)—on the system's evolution, providing nuanced insights into the cognitive and social foundations of online collective attention.
\end{enumerate}

The remainder of this paper is structured as follows. Section 2 reviews related work in topic modeling, agent-based simulation, and evolutionary game theory, positioning our integrative approach within the broader literature. Section 3 details our methodological pipeline, including data processing, topic extraction with BERTopic, payoff matrix construction, and the agent-based evolutionary game model. Section 4 presents our simulation experiments and evaluation on two Weibo datasets, analyzing the effects of different parameters and payoff structures. Finally, Section 5 concludes by discussing the implications of our findings, acknowledging limitations, and outlining future research directions for enhancing the model's realism and predictive power in modeling complex social processes.

\section{Related Work}

The challenge of modeling the complex interplay between information diffusion and opinion evolution on social media requires a synthesis of paradigms that have traditionally been developed in isolation across computational linguistics, social psychology, and complex systems science. Our work is situated at this interdisciplinary nexus. This section provides a critical review of the foundational literature, first examining the behavioral modeling paradigms of Agent-Based Modeling and Evolutionary Game Theory, then tracing the evolution of content analysis through topic modeling, and finally, synthesizing previous integration attempts to precisely define the methodological and theoretical gap our framework is designed to fill.
\subsection{Modeling Social Behavior: From Rationality to Bounded Rationality and Social Learning}

A central pursuit in computational social science is understanding how micro-level individual actions and local interactions aggregate into emergent macro-level collective phenomena. Early models, while foundational, were often analytically rigid. Threshold models of collective behavior \cite{granovetter1978} provided a powerful intuition for social contagion but struggled with agent heterogeneity. Similarly, the two-step flow of communication theory \cite{katz1955} highlighted the role of opinion leaders but lacked the formal mechanics to simulate dynamic, multi-step processes in digital networks \cite{lerman2010}.
\subsubsection{Agent-Based Modeling: A Generative Paradigm for Complex Social Systems}

Agent-Based Modeling (ABM) represents a paradigm shift towards a generative approach to social science \cite{epstein1996}, enabling the study of emergent phenomena from the bottom up \cite{helbing2012agent}. By simulating the actions and interactions of autonomous, heterogeneous agents, ABM captures the inherent complexity of social media, characterized by local rules, nonlinear dynamics, and decentralized control. Its application to online platforms has been extensive, from simulating information cascades on Digg and Twitter \cite{lerman2010} and modeling behavioral contagion in networked environments \cite{el2012social} to understanding the formation of polarization and echo chambers \cite{sasahara2021}. The flexibility of ABM is facilitated by platforms like NetLogo \cite{tisue2004netlogo}, which offers an accessible environment for modeling social learning, and specialized frameworks like SCRIMMAGE \cite{demarco2019simulating} and Cougaar \cite{helsinger2004cougaar} used in more complex, distributed simulations. Tools like the Altreva Adaptive Modeler \cite{marica2015simulating} further demonstrate the application of ABM in simulating realistic market behaviors, underscoring the method's versatility. Recent efforts have enhanced ABM's realism for information diffusion by integrating semantic topic dynamics, such as modeling how user interest evolves and how derivative topics emerge from primary ones \cite{li2025dynamic}, or by incorporating network inference techniques to uncover latent connections that influence propagation \cite{jia2025rumor}. However, the strength of ABM—its flexibility—can also be a weakness; without a rigorous theoretical foundation for agent decision-rules, models risk being ad-hoc "story generators."
\subsubsection{Evolutionary Game Theory: Formalizing Adaptation and Bounded Rationality}

Evolutionary Game Theory (EGT) provides the necessary mathematical rigor to formalize agent decision-making, moving beyond the hyper-rationality of classical game theory \cite{kandori1997evolutionary, smith1973logic}. EGT models populations of boundedly rational agents who learn and adapt through mechanisms like imitation and social comparison, rather than utility optimization \cite{kandori1997evolutionary, gigerenzer1996}, \cite{hofbauer1998evolutionary}, \cite{weibull}. This theoretical shift is profoundly relevant to social media, where users often rely on simple heuristics, imitating behaviors that appear successful—a process directly aligned with social learning theory \cite{bandura1971} and the study of social learning strategies \cite{rendell2010}, \cite{sandholm2010population}. EGT formalizes this through revision protocols and pairwise comparison dynamics \cite{weibull}, \cite{hofbauer1998evolutionary}, \cite{sandholm2010population}, which dictate how agents switch strategies based on observed payoffs \cite{Broom2022}. This makes EGT exceptionally apt for modeling topic competition, where different themes or narratives can be treated as competing strategies, and their "fitness" is dynamically determined by user engagement metrics \cite{lorits2023evolutionary, wei2021online}. Applications of EGT are diverse \cite{Wang2017}, spanning the evolution of cooperation \cite{nowak2006}, social norms \cite{axelrod1986}, supply chain logistics \cite{zhiwen2020supply}, and public goods dilemmas \cite{kandori1997evolutionary}. EGT has also been successfully applied to model strategic interactions in rumor propagation, formally analyzing the roles of users, media, and government \cite{li2025rumor}. Furthermore, it informs models designed to guide information toward specific target audiences through adaptive incentive mechanisms \cite{meng2025modeling}. The integration of EGT with ABM creates a powerful synergy: ABM provides the "stage" for interactions, and EGT provides the psychologically plausible "script" for decision-making. Crucially, both are anchored in the concept of bounded rationality \cite{simon1990, kahneman2003}, which we operationalize through behavioral noise, simulating the cognitive constraints, imperfect judgment, and accelerated dynamics of collective attention \cite{lorenzspreen2019} inherent in real-user behavior. Pushing this integration further, recent models explicitly couple network game theory with opinion dynamics to create Public Opinion Management Models (POMM). These frameworks are designed  to simulate the co-evolution of individual behavior and opinion and proactively test the efficacy of public opinion management strategies \cite{xu2025pomm}.
\subsection{Understanding Content: The Semantic Revolution in Topic Modeling}

While behavioral models simulate the mechanisms of information spread, a complete understanding requires analyzing the semantic content that is being propagated. Topic modeling has been the cornerstone technique for uncovering latent thematic structures in large, unstructured text corpora \cite{egger2022topic, blei2012}.
\subsubsection{The Probabilistic and Matrix Factorization}

The field was fundamentally shaped by probabilistic generative models. Probabilistic Latent Semantic Analysis (PLSA) \cite{hofmann1999probabilistic} introduced the concept of documents as mixtures of latent topics. This was significantly advanced by Latent Dirichlet Allocation (LDA) \cite{blei2003latent}, which became the workhorse of the field by adding Dirichlet priors, providing a more complete generative model. Subsequent refinements introduced hierarchical structures (e.g., hLDA \cite{yu2019hierarchical}) and the Pachinko Allocation Model (PAM) \cite{li2006pachinko} to capture topic correlations. In parallel, matrix factorization approaches like Non-negative Matrix Factorization (NMF) \cite{arora2012learning} and Latent Semantic Analysis (LSA) \cite{dumais2004latent} offered computationally efficient alternatives. However, these traditional models, reliant on bag-of-words representations, often struggle with the short, noisy, and context-dependent nature of social media text, leading to topics that can be semantically incoherent or difficult to interpret \cite{egger2022topic, blekanov2020}.
\subsubsection{The Neural Revolution: Contextualized Embeddings and Modern Topic Models}

The advent of transformer-based Large Language Models (LLMs) like BERT \cite{koroteev2021bert} and RoBERTa \cite{liu2019roberta} marked a revolution in NLP, generating deep, contextualized word embeddings that capture semantic meaning with unprecedented accuracy. This breakthrough enabled a new generation of topic models that overcome the limitations of their predecessors.

BERTopic \cite{grootendorst2022} is a leading model in this new lineage. It is not a single algorithm but a flexible, modular pipeline that leverages pre-trained transformer models for document embedding, uses dimensionality reduction techniques, and applies density-based clustering (HDBSCAN) to identify topical clusters in a way that is highly robust to noise. Its use of a class-based variant of TF-IDF to create topic representations results in coherent and interpretable topics. This approach is particularly effective for social media data, as demonstrated in applications mapping news topics to audience engagement metrics \cite{jiang2024storytelling} and in tracing the dynamic evolution and bifurcation points of online discussions, effectively modeling opinion cumulation \cite{blekanov2023mapping}. The versatility of pipelines like BERTopic is also evidenced by its broad application across social network analysis as a whole for tasks ranging from fake news detection to public health analysis, as highlighted in systematic reviews of AI techniques in the domain \cite{nedungadi2025ai}. The semantic coherence of BERTopic's output is critical for our framework, as it allows us to move beyond treating topics as mere statistical artifacts and instead define them as meaningful "strategic options" or "narratives" that agents can adopt and propagate, thereby creating a robust bridge between the semantic and behavioral layers of online discourse. It is also important to note that, while not strictly necessary for this research, robust evaluation techniques are lacking and existing approaches fail to fully assess a model's ability to induce a meaningful semantic structure on the document space. This has been explored as another area of combining game theory and text modeling as novel frameworks propose addressing this issue by integrating clustering evaluation metrics and unifying multi-perspective assessments through game-theoretic approaches like Multi-Attribute Utility Theory (MAUT) \cite{pereira2025current}.

\subsection{The Integration Gap: A Critical Synthesis of Previous Attempts and Remaining Challenges}

Recognizing the limitations of isolated approaches, a growing body of literature has attempted to integrate content analysis with behavioral simulation. A systematic review, however, reveals that these attempts have often resulted in a superficial coupling, failing to capture the dynamic, co-evolutionary relationship between semantic meaning and social action.

Early and some contemporary efforts involve using the output of NLP models as static inputs or fixed parameters for behavioral models. For instance, studies have used LDA to pre-define topics which then act as immutable categories of information that agents can share or be influenced by \cite{almaatouq2021}. Other works have incorporated sentiment analysis into ABMs to simulate emotion contagion \cite{fan2018, wu2025, mir2025joint, chen2025computational} or used network and content features to predict information cascade sizes \cite{cheng2014}. While these works provide valuable insights, they suffer from three interconnected limitations that our framework is designed to overcome:

\begin{enumerate}
\item Weak Semantic-Behavioral Integration and Co-evolution. In many models, the content layer (topics/sentiments) is a static label. It does not co-evolve with the agent population. The topics themselves do not possess dynamic "fitness" that changes based on collective agent behavior, and the mechanisms that drive the rise and fall of topics are not explicitly modeled within a unified simulation framework. This creates a critical disconnect where content may influence behavior, but behavior does not recursively shape the semantic landscape \cite{cheng2014}.

\item Use of Outdated or Semantically Shallow NLP Methods. Applying traditional models like LDA to short-text social media platforms often yields topics that are non-interpretable or fail to accurately reflect the nuanced, contextual discourses present in the data \cite{egger2022topic, grootendorst2022}. This "garbage in, garbage out" problem fundamentally undermines the validity of the subsequent behavioral simulation. If the "strategies" available to agents are poor representations of real-world discourse, the model's explanatory power is severely limited.

\item Lack of Psychological and Theoretical Plausibility. Agent decision-rules are frequently based on simple, mathematically convenient heuristics that lack grounding in established psychological theories. There is often a neglect of social learning \cite{bandura1971, rendell2010}, bounded rationality \cite{simon1990, gigerenzer1996}, and the cognitive constraints of information processing \cite{kahneman2003}. The concept of "cumulative deliberation" \cite{bodrunova2019, blekanov2023mapping}—which describes how public opinion online is formed through the accumulation of often non-deliberative individual expressions rather than rational debate—is rarely operationalized in computational models. Furthermore, the role of behavioral noise, imperfect imitation, and the exploration-exploitation trade-off as fundamental drivers of diversity and dynamism in real systems is often overlooked or treated as a mere calibration parameter rather than a core component of the theoretical framework \cite{lorits2023evolutionary}.
\end{enumerate}

Our proposed framework is designed to bridge this critical gap. We perform a deep integration of state-of-the-art, LLM-powered topic modeling (BERTopic) with a psychologically plausible ABM/EGT framework. This is not a loose coupling but a unified approach: the semantically coherent topics extracted by BERTopic become the dynamic strategies in an evolutionary game. Agents update their topic preferences through mechanisms of social learning (imitation) under bounded rationality (noise), directly reflecting the cognitive and social processes of users in online environments. This creates a closed, generative loop where the semantic content of topics influences individual agent behavior, and the aggregate behavior of the population, in turn, shapes the evolutionary trajectory and competitive success of the topics themselves, offering a more holistic and psychologically grounded model of online discourse.

\section{Method} \label{Method}

\subsection{Pipeline Architecture}
\begin{figure}
    \centering
    \includegraphics[width=.8\linewidth]{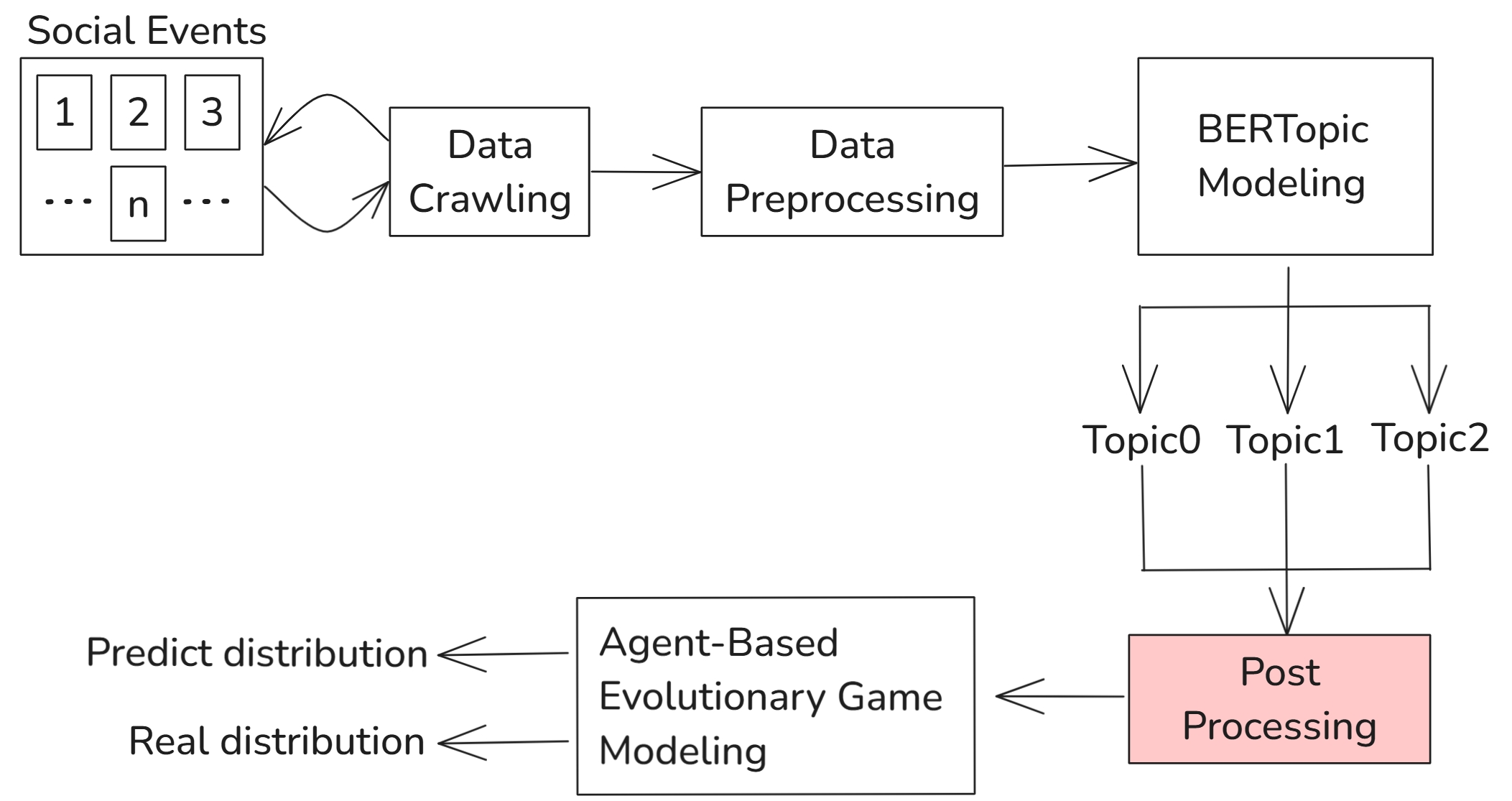}
    \caption{Pipeline}
    \label{fig:pipeline}
\end{figure}

The pipeline of this study is shown in Figure \ref{fig:pipeline}. First, data is collected from the Weibo platform using web crawling techniques \cite{khder2021web}, and the data is preprocessed. Subsequently, the preprocessed data is fed into the BERTopic model to extract topics and their related features. Then, the process moves to Post Processing, where, after a series of treatments, it provides the input for the subsequent multi-agent simulation. Next, a multi-agent model is built using NetLogo to simulate the dynamic dissemination process of topics within the social network. Finally, the original and predicted data distributions are obtained and they are analysed and compared.

\subsection{Data Processing and Topic Extraction}
Weibo posts from October and December 2024 were collected in six-hour polling windows and deduplicated via the triplet {content, user ID, timestamp} to yield a redundancy-free corpus; after character cleaning and time normalization, the corpus served as input to BERTopic and the influence metric.

Weibo posts were embedded with bert-base-chinese, dimension-reduced by PCA, and clustered with HDBSCAN; salient keywords were extracted via C-TF-IDF and manually labelled to define agent strategies.

\subsection{Post Processing Module}
This module is one of the core components of the process, and its primary task is to transform the topic information extracted by BERTopic into input data suitable for multi-agent simulation. The specific steps are as follows:
\begin{itemize}
\item[1.] Information influence is quantified by ranking each post based on descending order of likes, comments, and reposts. Due to the typically skewed distribution of engagement metrics (where most content receives minimal interactions while a few viral posts achieve exceptionally high numbers), conventional equal-frequency binning methods prove inadequate as they create artificial groupings. We instead employ a manual stratification approach based on the intrinsic characteristics of the data distribution. This involves first analyzing key statistical features of each metric (including extreme values, means, and percentiles) to identify natural breakpoints in the data. These insights then inform the division of content into 10 distinct tiers with meaningful engagement thresholds.

The methodology offers three key advantages: it preserves authentic distribution patterns, yields more interpretable categorizations for practical applications, and effectively captures outlier viral content. This approach produces rankings that more accurately reflect real-world influence dynamics compared to algorithmic quantization methods.


\item[2.] The payoff or utility of a single piece of information is calculated using the next steps. 


For each post, the numbers of likes, comments, and reposts are separately ranked from highest to lowest, and the ranking intervals are divided into ten decile levels (1 = lowest, 10 = highest). Let $ 
R_l,\; R_c,\; R_r \in \{1,2,\dots,10\}
$
denote these decile ranks for a given post. The individual post payoff is defined as the average rank:
\begin{equation}\label{post_payoff}
q = ( R_l + R_c + R_r)/3.
\end{equation}
This value lies in \([1,10]\) and represents the post's overall engagement influence. Next, for each topic \(i \in \{1,2,3\}\) that contains \(N_i\) posts, we compute its average topic payoff $p_i$ as
\begin{equation}\label{topic_payoff}
    p_i = \frac{1}{N_i} \sum_{\text{posts } k \in \text{Topic}_i} q_k,
\end{equation}
where \(q_k\) is the individual payoff of post \(k\). These \(p_i\) quantify the aggregate influence of each topic and serve as the base payoffs for the evolutionary game.




To further explore the relative contributions of different interaction dimensions, we introduce weighting parameters $\beta, \gamma, \delta \in [1,3]$ for likes, comments, and reposts, respectively. For each post, its weighted individual payoff is defined as
\begin{equation}
\tilde{q} = \frac{\beta R_l + \gamma R_c + \delta R_r}{3},
\label{eq:payoff_weighted}
\end{equation}
where $R_l, R_c, R_r \in \{1,\dots,10\}$ are the decile ranks of that post. Correspondingly, for topic $i$ containing $N_i$ posts, the weighted average topic payoff is given by
\begin{equation}\label{eq:weighted_topic_payoff}
\tilde{p}_i = \frac{1}{N_i} \sum_{k \in \text{Topic}_i} \tilde{q}_k,
\end{equation}
where $\tilde{q}_k$ is the weighted individual payoff of post $k$ computed according to \eqref{eq:payoff_weighted} with the corresponding ranks. By varying these weights, we can simulate diverse social contexts where different engagement metrics carry distinct strategic importance. This extended formulation provides a quantitative basis for assessing the relative effects of likes, comments, and reposts on information diffusion within the evolutionary game framework, while directly generalising the base case $q_k$ and $p_i$ (which correspond to $\beta=\gamma=\delta=1$).

The topic payoffs \(p_i\) (or \(\tilde{p}_i\) for the weighted case) obtained from \eqref{topic_payoff} (or \eqref{eq:weighted_topic_payoff}) are then used as the base influence values \(P_i\) for constructing the payoff matrix in the agent-based simulation (see Section~3.5).




\end{itemize}

In summary, the Post Processing Module realizes the standardized transformation from raw textual data to simulation-ready inputs. It not only assigns a quantitative influence characterization to each topic but also, through the adjustable parameter mechanism, lays a solid foundation for sensitivity analyses and strategic exploration in subsequent simulation experiments.

\subsection{Agent-Based Evolutionary Game Modeling}

We analyze an evolutionary game in which a finite population of agents engages in repeated random pairwise interactions. Each agent is programmed with a pure strategy, corresponding to a discrete action (specifically, following a particular topic). The payoff from these symmetric encounters determines the strategy's fitness. The population state evolves via a strategy revision process: in discrete time steps, agents may switch strategies by comparing their realized payoffs, preferentially imitating more successful behaviors. This dynamic governs the evolutionary trajectory of strategy frequencies within the population.

The population distribution is characterized by strategy frequencies. For each strategy (topic) \(i\), the population share \(x_i(t)\) at time \(t\) is given by:
\begin{equation}
x_i(t) = \frac{n_i(t)}{N},
\label{eq:population_share}
\end{equation}
where \(n_i(t)\) denotes the number of agents following strategy \(i\) at time \(t\), and \(N\) is the total population size. 

The strategy set \(\mathcal{S} = \{0, 1, \dots, n\}\) corresponds to the topics identified via BERTopic clustering. The initial population state is defined by a vector \(\mathbf{x}(0) = (n_0(0), n_1(0), \dots, n_n(0))\), where \(\sum_i n_i(0) = N\).

The game is characterized by a symmetric payoff matrix that specifies the fitness consequences for both participants in an encounter. Let \( \mathbf{A} = (A_{ij})_{i,j \in \mathcal{S}} \) be the payoff matrix. When an agent using strategy \( i \) interacts with an agent using strategy \( j \), the former receives payoff \( A_{ij} \), while the latter receives payoff \( A_{ji} \). The matrix is typically not symmetric (\( A_{ij} \neq A_{ji} \) in general).

\[
\mathbf{A} = \begin{bmatrix}
A_{00} & A_{01} & \cdots & A_{0n} \\
A_{10} & A_{11} & \cdots & A_{1n} \\
\vdots & \vdots & \ddots & \vdots \\
A_{n0} & A_{n1} & \cdots & A_{nn}
\end{bmatrix}.
\]

The dynamics are driven by a \textbf{pairwise imitation protocol with noise} \citep{sandholm2010population}. In each discrete time period \(t\), the following steps occur:

\begin{enumerate}
    \item \textbf{Pairwise Interaction:} Every agent is randomly paired with another agent from the population. Consider an agent using strategy \(i\) paired with an agent using strategy \(j\). According to the payoff matrix \(\mathbf{A}\), the first agent receives payoff \(A_{ij}\), while the second receives \(A_{ji}\). Let \(\pi_i = A_{ij}\) denote the payoff received by the agent with strategy \(i\).
    
    \item \textbf{Strategy Revision:} With probability \(\rho \in [0,1]\) (the \textit{revision rate}), each agent receives an opportunity to revise its strategy. Consider a revising agent that used strategy \(i\) and received payoff \(\pi_i\). This agent selects a random opponent from the population who used strategy \(j\) and received payoff \(\pi_{j} = A_{jk}\) for some \(k \in \mathcal{S}\).
    \begin{itemize}
        \item With probability \(\eta \in [0,1]\) (the \textit{noise level}), the agent adopts a strategy uniformly at random from \(\mathcal{S}\).
        \item With probability \(1-\eta\), the agent imitates the opponent's strategy \(j\) if and only if \(\pi_{j} > \pi_i\) (the ``imitate-if-better'' rule).
    \end{itemize}
\end{enumerate}

All strategy updates are applied synchronously at the end of the time step. This process defines an evolutionary dynamic on the space of population states.

The Algorithm~\ref{alg:update} provides a detailed pseudocode implementation of one time step (tick) of this evolutionary update process.

\begin{algorithm}[H]
\caption{One Tick of the Evolutionary Update Process}
\label{alg:update}
\begin{algorithmic}[1]
\State \textbf{Input:} Agent population $\mathcal{A}$, payoff matrix $A$, revision probability \texttt{prob-revision}, noise level \texttt{noise}

\For{each agent $a_i \in \mathcal{A}$}
    \State Randomly select another agent $a_j \in \mathcal{A}, j \neq i$
    \State Agent $a_i$ plays a game with $a_j$ and obtains payoff
\EndFor

\For{each agent $a_i \in \mathcal{A}$}
    \If{revision occurs with probability \texttt{prob-revision}}
        \If{random() $<$ \texttt{noise}}
            \State $a_i$ adopts a random strategy
        \Else
            \If{payoff$(a_j) >$ payoff$(a_i)$}
                \State $a_i$ imitates the strategy of $a_j$
            \EndIf
        \EndIf
    \EndIf
\EndFor

\State Apply all strategy updates \textbf{synchronously} at the end of the tick

\end{algorithmic}
\end{algorithm}

\subsection{Parameter Specification}

The simulation uses three tunable inputs.

\begin{enumerate}
  \item \textbf{Payoff matrix.}
    \begin{itemize}
      \item Additive method: $A_{ij}=P_i+P_j$.
      \item Weighted method: $A_{ij}=\alpha P_i+(1-\alpha)P_j$.
    \end{itemize}
    where $P_i$ denotes the average influence of topic $i$, \( \alpha \in [0,1] \) adjusts the contribution ratio of self-strategy (\( P_i \)) and other-strategy (\( P_j \)) to the agent's payoff. The payoff matrix, which quantifies the strategic interactions between agents, explicitly shows the relative importance of different topics in the evolutionary game, as it is proposed in \eqref{topic_payoff}-\eqref{eq:weighted_topic_payoff}. Specifically, the matrix elements represent the fitness payoffs when an agent following topic $i$ interacts with an agent following topic $j$.

  \item \textbf{Revision probability.}
    It denotes the probability that an agent attempts to revise its current strategy in each simulation round, and can be interpreted as the "propensity for change" or "information update frequency" of social media users. Extensive simulation results reveal that lower values of prob-revision allow the model to better fit the observed temporal trends in real topic proportions. Based on this observation, comparative experiments systematically test prob-revision values of 0.03, 0.06, and 0.12 to evaluate differences in model stability and fitting performance under low-frequency (slow evolution) and high-frequency (rapid adjustment) change scenarios.

  \item \textbf{Behavioural noise.}
    It is introduced to characterize the irrational perturbations in agent decision-making when considering strategy revisions.  It helps prevent the system from becoming trapped in local optima and enhances the diversity and dynamism of evolutionary dynamics. Theoretically, a higher noise value increases system diversity and delays early convergence trends, whereas a lower noise value accelerates strategy convergence and mainstream adoption within the population. Empirical testing indicates that the overall influence of the noise parameter on the trend is relatively limited; therefore, in most simulation experiments, noise is uniformly set to 0.4 as a reasonable approximation to model bounded rationality and information misjudgment in social environments.
\end{enumerate}

Parameter sets are chosen by minimising RMSE and nRMSE \cite{yang2014evaluation} between simulated and observed topic shares.

\section{Simulation and evaluation experiments} \label{Experiments}

\subsection{Experimental Setup}

The dataset selected in this study centers on the topic of \textit{climate change}, based on the following three main considerations. First, from the perspective of temporal distribution of public attention, climate-related discussions occur more frequently during the autumn and winter seasons. Second, the topic of climate change encompasses a rich semantic structure, covering multiple subtopics. Third, climate change is both highly realistic and globally relevant.

Based on this theme, two real-world social media datasets from the Weibo platform were constructed in this study, serving respectively to validate the feasibility of the model and to evaluate its generalization ability. The first dataset covers the period from October 6, 2024, to October 31, 2024, and contains approximately 20,000 Weibo posts. The second dataset covers December 6, 2024, to December 31, 2024, comprising around 12,000 posts.

The NetLogo platform was used to simulate individual behaviors and topic evolution processes. For Dataset 1, the model is initialized with 310 agents, with a strategy distribution of  \([249, 34, 27]\), corresponding to the number of posts(agents) assigned to three topics by BERTopic within the initial time window. For Dataset 2, the initialization involved 340 agents with a strategy distribution of \([250, 58, 32]\). In both datasets, the noise parameter was uniformly set to 0.4 to simulate uncertainty in real-world social behavior. Dynamic variables, such as the strategy revision probability and payoff matrix structures, will be further adjusted and compared in subsequent sections.

\subsection{Case Study 1: October 2024 Climate Change Dataset}

\subsubsection{Topic Structure}

After applying BERTopic modeling to the Weibo corpus collected from October 6, 2024, to October 31, 2024, three representative topics were extracted, which were subsequently named ``temperature changes'', ``air pollution'', and ``extreme weather.'' It should be noted that the topics output by BERTopic under the default configuration are not presented as complete natural language expressions but rather as clusters of high-frequency keywords (e.g., ``0\_Temperature \_Today\_Daylight\_Area''). To enhance the clarity and readability of the subsequent analysis, the preliminary outputs from BERTopic were further interpreted. Specifically, based on the semantic content of the top keywords associated with each topic, this study manually refined and consolidated the topics into semantically representative short-phrase names.

Table~\ref{tab:subtopic_summary} summarizes the classification results obtained using the BERTopic model on the Weibo dataset, as well as the corresponding payoff values for each topic.

\begin{table}[ht]
\centering
\caption{Topic Classification Results and Payoff}
\label{tab:subtopic_summary}
\begin{tabular}{|c|c|c|c|}
\hline
\textbf{Topic} & \textbf{Keyword} & \textbf{Count} & \textbf{Payoff} \\
\hline
0 & Temperature changes & 12413 & 2.085529 \\
\hline
1 & Air pollution & 4191 & 1.518253 \\
\hline
2 & Extreme weather & 3952 & 2.110240 \\
\hline
\end{tabular}
\end{table}

\begin{figure}
    \centering
    \includegraphics[width=.9\linewidth]{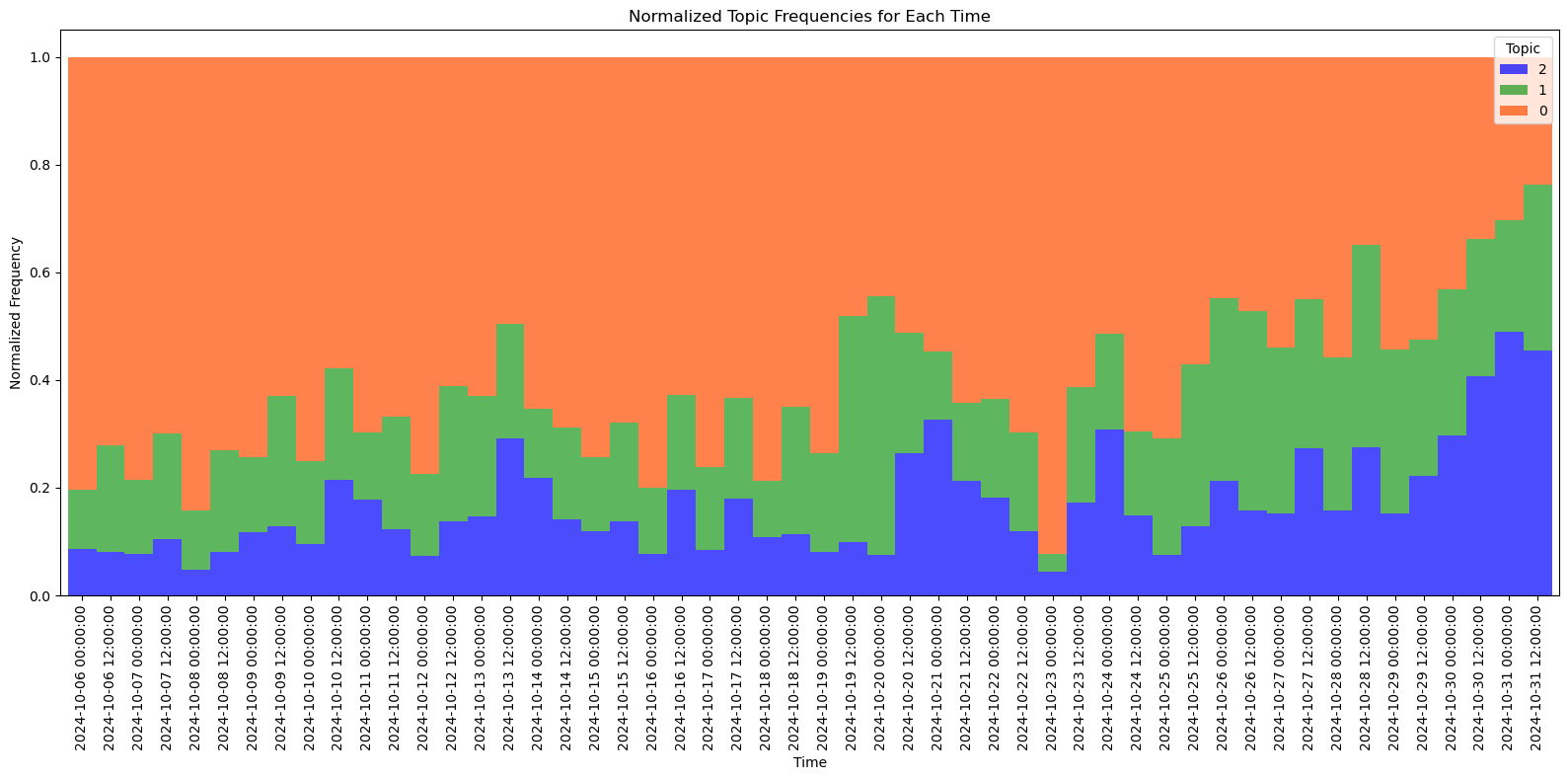}
    \caption{Distribution of raw data}
    \label{fig:Distribution of raw data}
\end{figure}

 To further analyze the topic evolution trends, Figure~\ref{fig:Distribution of raw data} presents the distribution of the three topics sampled every 12 hours from the original data, covering a total of 52 time points.

However, due to the inherent presence of noise and outliers in social media data, such as sudden bursts of reposts or abnormal spikes in comments, directly using the raw data for modeling could introduce significant fluctuations, adversely affecting simulation stability and fitting accuracy. Therefore, the original time series data were smoothed prior to modeling. Specifically, a Moving Average (MA) technique \cite{hansun2013new} was applied to reduce randomness and highlight long-term trend variations. To balance the trade-off between smoothing strength and sensitivity to local changes, this study explored three commonly used window sizes: 4, 5, and 7.

\begin{figure}
    \centering
    \includegraphics[width=.9\linewidth]{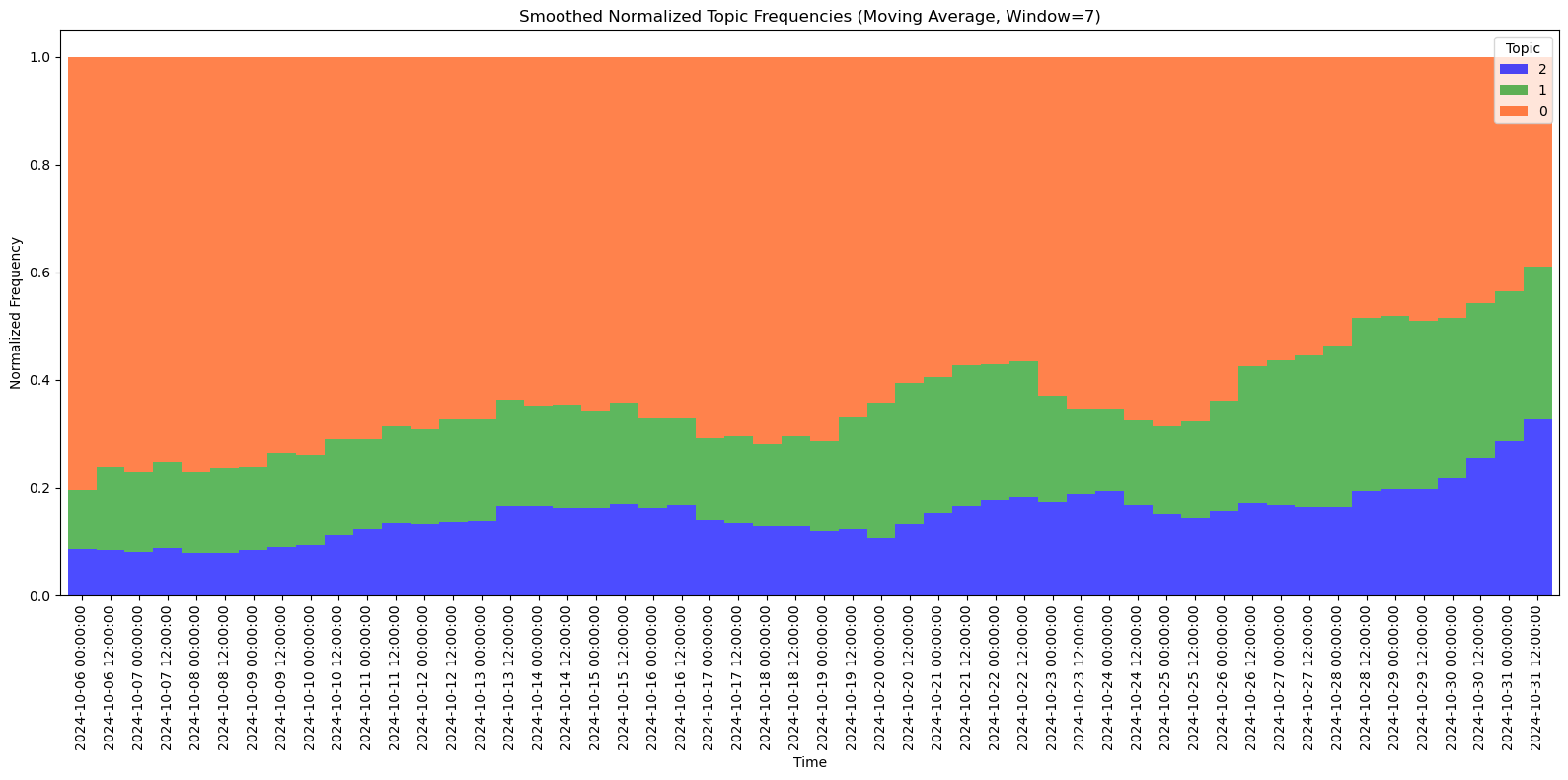}
    \caption{Distribution of raw data after moving average smoothing (window size 7)}
    \label{fig:Window_size_7}
\end{figure}

Figure~\ref{fig:Window_size_7} presents the distribution of the three topics after applying Moving Average smoothing with a window size of 7. It can be observed that after smoothing, short-term fluctuations are significantly reduced, and the trend variations become more continuous and clearer, thereby providing a stable foundation for subsequent NetLogo behavioral simulations and comparisons with real-world data.

\subsubsection{Influence of Prob-revision Value on the Fitting Effect}

In this set of experiments, we evaluate the impact of the individual strategy revision probability (\texttt{prob-revision}) on the dynamic evolution of the model. Under a fixed payoff matrix, simulations were conducted with three different \texttt{prob-revision} values: 0.03, 0.06, and 0.12, respectively, to perform comparative analysis. The payoff matrix used for this experiment was constructed using the additive method, which is designated as the baseline configuration. The resulting payoff matrix is as follows:

\begin{equation}
A =
\begin{bmatrix}
4.171058 & 3.603782 & 4.195769 \\
3.603782 & 3.036506 & 3.628493 \\
4.195769 & 3.628493 & 4.220480
\end{bmatrix}
\label{eq:baseline_payoff_matrix}.
\end{equation}

Under this fixed payoff structure, each agent in the social network observes the strategy performance of neighboring agents with a certain probability and updates its own strategy by imitation if a higher expected payoff is observed. By adjusting the \texttt{prob-revision} parameter, the degree of agents' ``learning'' or ``updating'' activity can be controlled, thereby influencing the overall evolution trajectory of strategies.

The experimental results are summarized in Table~\ref{tab:prob_revision_comparison}. When \texttt{prob-revision} was set to 0.06, the model achieved the lowest error between the simulated and real topic evolution paths, with an RMSE of 0.114 and an nRMSE of 12.8\%. In comparison, when \texttt{prob-revision} was 0.03, the model exhibited a more conservative behavior with lower imitation frequency. Although initial stability was stronger, the model was slightly slower in responding to strategy evolution, leading to a modest increase in error (RMSE = 0.122, nRMSE = 13.7\%). When \texttt{prob-revision} was increased to 0.12, frequent individual updates caused higher system volatility and instability, resulting in greater deviations from real trends, with RMSE rising to 0.188 and nRMSE reaching 21.1\%.

\begin{table}[ht]
\centering
\caption{Comparison of Simulation Errors under Different \texttt{prob-revision} Settings}
\label{tab:prob_revision_comparison}
\begin{tabular}{|c|c|c|}
\hline
\textbf{\texttt{prob-revision}} & \textbf{RMSE} & \textbf{nRMSE} \\
\hline
0.03 & 0.122 & 13.7\% \\
\hline
0.06 & \textcolor{red}{0.114} &\textcolor{red}{12.8\%} \\
\hline
0.12 & 0.188 & 21.1\% \\
\hline
\end{tabular}
\end{table}

In further experiments, the errors between the model simulation results and the real-world topic evolution paths smoothed by the Moving Average (MA) method were computed. The experimental results are summarized in Table~\ref{tab:smoothing_comparison}. When using a window size of 7 and setting \texttt{prob-revision} to 0.03, the model achieved the smallest error, with an RMSE of 0.085 and an nRMSE of 11.7\%, which clearly outperformed other combinations of window sizes and revision probabilities. This finding suggests that a window size of 7 strikes an optimal balance between suppressing short-term fluctuations and preserving the true evolutionary rhythms. Given that a natural weekly cycle consists of 7 days, using 7 as the smoothing window size may better align with potential periodic patterns inherent in social topic evolution. Therefore, unless otherwise specified, a smoothing window size of 7 is adopted as the standard reference in subsequent comparisons and evaluations.

\begin{table}[ht]
\centering
\caption{Comparison of Simulation Errors under Different Moving Average Window Sizes}
\label{tab:smoothing_comparison}
\begin{tabular}{|c|cc|cc|cc|}
\hline
\multirow{2}{*}{\textbf{\texttt{prob-revision}}} & \multicolumn{2}{c|}{\textbf{Window Size = 4}} & \multicolumn{2}{c|}{\textbf{Window Size = 5}} & \multicolumn{2}{c|}{\textbf{Window Size = 7}} \\
\cline{2-7}
 & RMSE & nRMSE & RMSE & nRMSE & RMSE & nRMSE \\
\hline
0.03 & 0.095 & 13.1\% & 0.091 & 12.5\% & \textcolor{red}{0.085} & \textcolor{red}{11.7\%} \\
\hline
0.06 & 0.095 & 13.1\% & 0.094 & 12.9\% & 0.092 & 12.7\% \\
\hline
0.12 & 0.183 & 25.3\% & 0.184 & 25.4\% & 0.186 & 25.7\% \\
\hline
\end{tabular}
\end{table}

Overall, lower strategy revision probabilities (\texttt{prob-revision} = 0.06 or 0.03) achieved a better balance between system stability and dynamic responsiveness. Based on the preceding experimental results, the following settings were adopted for subsequent simulations on Dataset 1: when directly comparing the simulation outputs with the original raw data, \texttt{prob-revision} was set to 0.06 (corresponding to the optimal revision probability under the original data); when comparing the simulation outputs with the real data smoothed using a Moving Average window size of 7, \texttt{prob-revision} was set to 0.03 (corresponding to the optimal revision probability under the smoothed data). This configuration strategy ensures that the model achieves optimal fitting performance and explanatory power under different data processing scales, providing a unified and rational parameter foundation for the subsequent comparative experiments.

\subsubsection{Comparison of Payoff Matrix Methods}

After determining the optimal behavior revision probability in the previous section, this section further investigates the impact of different payoff matrix construction methods on model fitting performance.

To ensure fairness and comparability, other model parameters (such as \texttt{noise} = 0.4 and the initial strategy distribution \([249, 34, 27]\)) were kept consistent across simulations. For each experimental setting, simulations were conducted with \( \alpha \in \{0.1, 0.2, \ldots, 0.9\} \). The fitting performance under different payoff matrix structures was quantitatively evaluated by calculating the RMSE and nRMSE between the simulation outputs and the real-world data.

The experimental results for the first setting are summarized in Table~\ref{tab:alpha_original_results}. As shown in the results, when \( \alpha = 0.3 \), the model achieved the lowest error (RMSE = 0.098, nRMSE = 11.0\%), indicating that incorporating attention to the opponent’s strategy payoff in the payoff construction enhances the model’s ability to realistically capture information diffusion behavior in social networks, The simulation results are shown in  Figure~\ref{Simulation with a=0.3, prob_revision=0.06}.

Conversely, as \( \alpha \) gradually increased (i.e., as agents placed more emphasis on their own strategies), the model error consistently rose. Particularly when \( \alpha = 0.6 \), RMSE increased to 0.159, and the simulated trend significantly deviated from the real evolution. This suggests that insufficient responsiveness to others’ strategies leads to an overly ``self-reinforcing'' evolution, making it difficult to replicate the imitation and competition mechanisms observed in real systems. The baseline method (additive payoff construction without weighting) exhibited moderate overall performance (RMSE = 0.114), suggesting that although simple symmetric constructions have a certain degree of fitting capability, they are relatively limited in capturing dynamic differences among subtopics.

\begin{table}[ht]
\centering
\caption{Simulation Errors under Different \(\alpha\) Values (Original Data)}
\label{tab:alpha_original_results}
\begin{tabular}{|c|c|c|}
\hline
\textbf{\(\alpha\)} & \textbf{RMSE} & \textbf{nRMSE} \\
\hline
Baseline & \textcolor{orange}{0.114} & \textcolor{orange}{12.8\%} \\
\hline
0.1 & 0.104 & 11.6\% \\
\hline
0.2 & 0.107 & 12.0\% \\
\hline
0.3 &  \textcolor{red}{0.098} &  \textcolor{red}{11.0\%} \\
\hline
0.4 & 0.104 & 11.7\% \\
\hline
0.5 & 0.116 & 13.0\% \\
\hline
0.6 & \textcolor{blue}{0.159} & \textcolor{blue}{17.9\%} \\
\hline
0.7 & 0.129 & 14.4\% \\
\hline
0.8 & 0.140 & 15.7\% \\
\hline
0.9 & 0.120 & 13.5\% \\
\hline
\end{tabular}
\end{table}

\begin{figure}
    \centering
    \includegraphics[width=.9\linewidth]{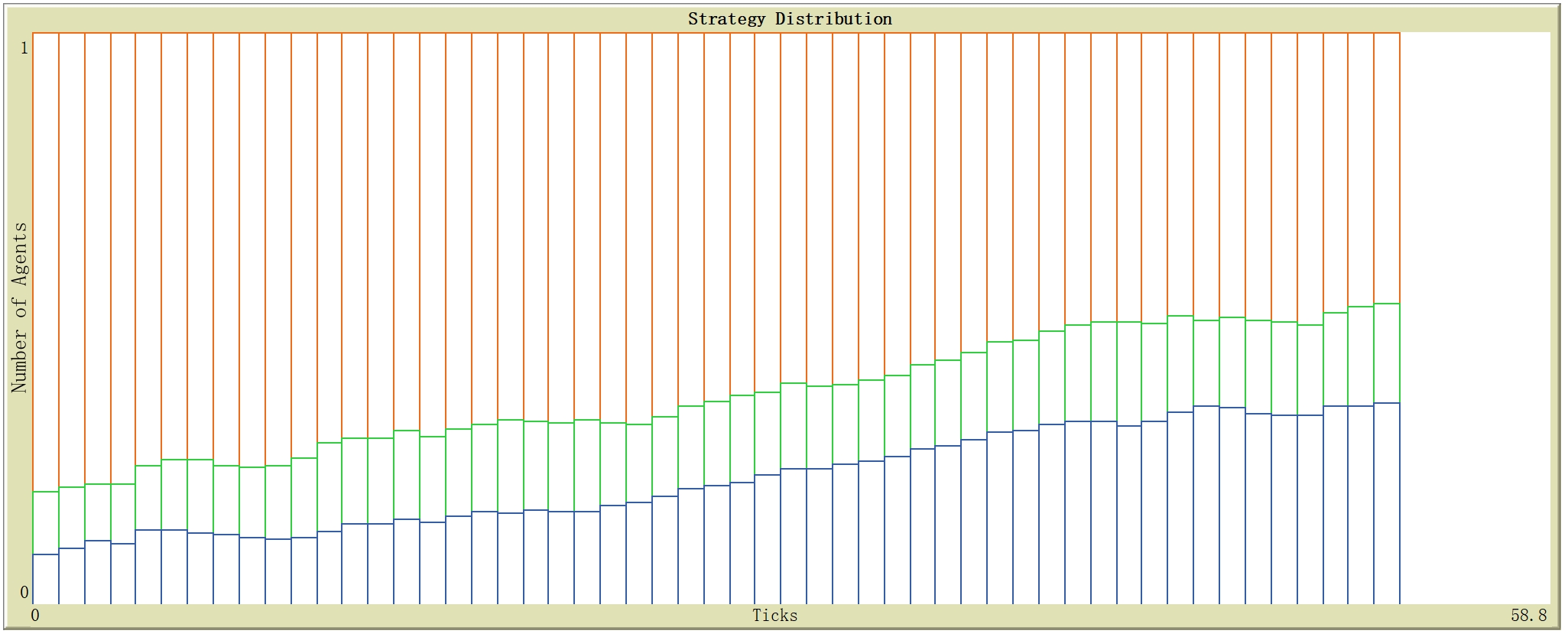}
    \caption{Simulation with \(\alpha\)=0.3, prob\_revision=0.06}
    \label{Simulation with a=0.3, prob_revision=0.06}
\end{figure}

In the setting with a smoothing window size of 7 and \texttt{prob-revision} fixed at 0.03, the simulation fitting results for different \( \alpha \) values are summarized in Table~\ref{tab:alpha_smoothed_results}. Compared to the previous experiments using raw (unsmoothed) data, the overall RMSE and nRMSE are significantly lower, indicating that the smoothing process effectively reduced noise fluctuations in the original data, thereby making the model fitting evaluation more stable and reliable. The simulated data distribution under this optimal setting is shown in Figure~\ref{Simulation with a=0.4, prob_revision=0.03}.

Moreover, from the overall trend, it can be observed that RMSE and nRMSE remained relatively low across the range \( \alpha \in [0.1, 0.9] \), reflecting the model’s robustness in the parameter space. Although slight variations in model performance were observed as \( \alpha \) changed, the best fitting consistently appeared around \( \alpha = 0.35 \). This phenomenon further validates the advantage of the weighted payoff structure in modeling individual decision trade-offs: it captures the dual driving forces of self-interest and responsiveness to environmental changes.

\begin{table}[ht]
\centering
\caption{Simulation Errors under Different \(\alpha\) Values (Smoothed Data)}
\label{tab:alpha_smoothed_results}
\begin{tabular}{|c|c|c|}
\hline
\textbf{\(\alpha\)} & \textbf{RMSE} & \textbf{nRMSE} \\
\hline
Baseline & \textcolor{orange}{0.085} & \textcolor{orange}{11.7\%} \\
\hline
0.1 & 0.082 & 11.4\% \\
\hline
0.2 & 0.074 & 10.3\% \\
\hline
0.3 & 0.076 & 10.6\% \\
\hline
0.4 & \textcolor{red}{0.055} & \textcolor{red}{7.6\%} \\
\hline
0.5 & 0.079 & 10.9\% \\
\hline
0.6 & 0.071 & 9.7\% \\
\hline
0.7 & 0.081 & 11.2\% \\
\hline
0.8 & 0.075 & 10.4\% \\
\hline
0.9 & 0.074 & 10.2\% \\
\hline
\end{tabular}
\end{table}

\begin{figure}
    \centering
    \includegraphics[width=.9\linewidth]{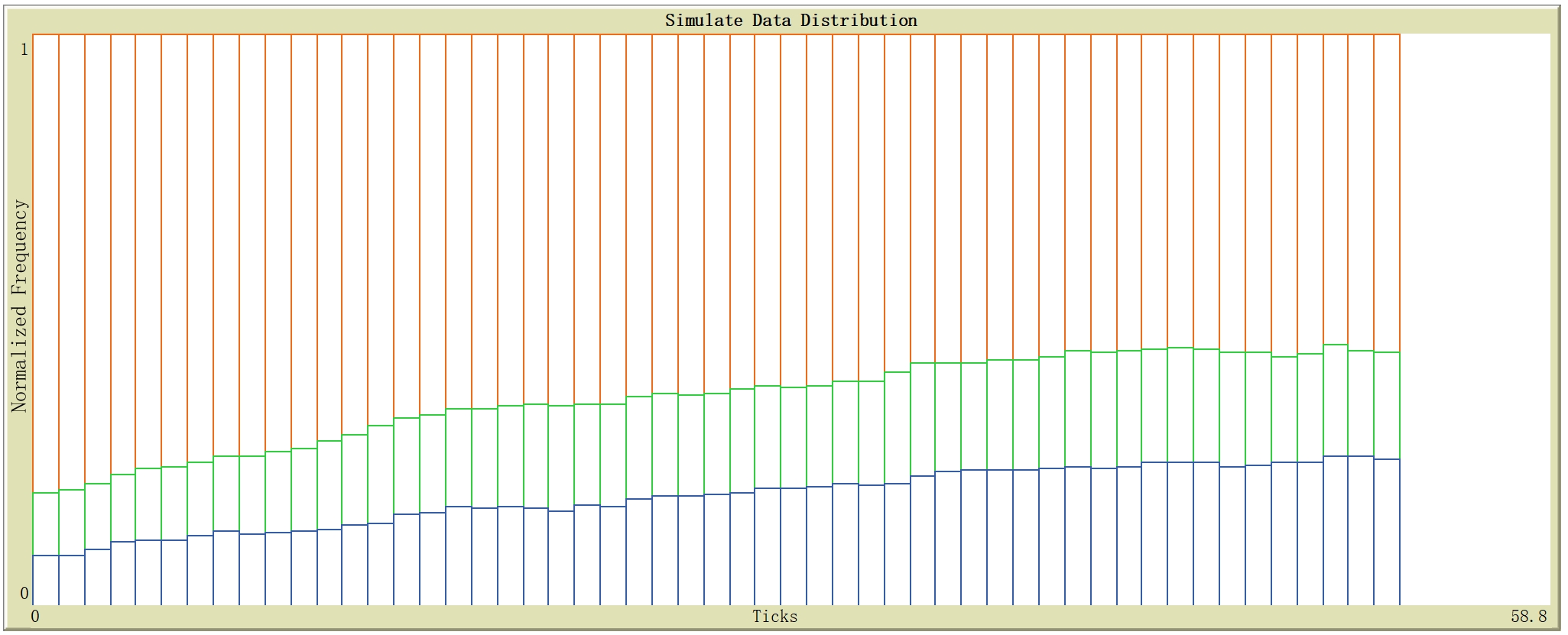}
    \caption{Simulation with \(\alpha\)=0.4, prob\_revision=0.03}
    \label{Simulation with a=0.4, prob_revision=0.03}
\end{figure}

Overall, the experiments conducted in Case Study 1 have fully validated the effectiveness and flexibility of the proposed integrated ``semantic modeling – evolutionary simulation'' framework. Through the coordinated optimization of multiple submodules, the model not only successfully reproduced the macro-level topic proportion changes but also provided explanatory insights into the micro-level behavioral evolution mechanisms of individual agents. To further assess the robustness and cross-dataset adaptability of the proposed framework, the next chapter will apply the same experimental procedures and evaluation methods to another Weibo dataset collected during a different time period (December 2024) focusing on weather variations.

\subsection{Case Study 2: December 2024 Climate Change Dataset}
\subsubsection{Topic Structure}

In this section, Weibo data collected from December 6, 2024, to December 31, 2024 (approximately 12,000 posts in total), were processed using the same BERTopic parameter configuration as in Case Study 1. Three major topics were extracted, and their occurrence frequencies and average payoff values were calculated. The results are summarized in Table~\ref{tab:subtopic_summary_december}.

\begin{table}[ht]
\centering
\caption{Topic Classification Results and Payoff (December Dataset)}
\label{tab:subtopic_summary_december}
\begin{tabular}{|c|c|c|c|}
\hline
\textbf{Topic} & \textbf{Keyword} & \textbf{Count} & \textbf{Payoff} \\
\hline
0 & Temperature changes & 6111 & 2.341679 \\
\hline
1 & Air pollution & 3392 & 1.549528 \\
\hline
2 & Extreme weather & 1956 & 2.293286 \\
\hline
\end{tabular}
\end{table}

\subsubsection{Model Parameter Tuning}

In this experiment, the payoff matrix construction method was fixed to the first type (additive structure), serving as the baseline for exploring the impact of different revision probability values on simulation fitting performance. This setup aims to provide a stable foundation for subsequent comparisons of different payoff mechanisms.

As shown in Table~\ref{tab:prob_revision_december}, the revision probability was set to values from the set \{0.03, 0.06, 0.12\}. Under identical initialization conditions and noise level (\texttt{noise} = 0.4), the simulation model was run for each setting, and the corresponding RMSE and nRMSE values were recorded. It can be observed that when \texttt{prob-revision} was set to 0.06, the model achieved the best fitting performance across both metrics (RMSE = 0.120, nRMSE = 17.1\%), outperforming the settings with 0.03 and 0.12. This result suggests that, for the current dataset (December 2024), setting the strategy revision probability at a moderate level better captures the imitation and adjustment behavior of individuals regarding topic preferences.

\begin{table}[ht]
\centering
\caption{Simulation Errors under Different \texttt{prob-revision} Settings (December Dataset)}
\label{tab:prob_revision_december}
\begin{tabular}{|c|c|c|}
\hline
\textbf{\texttt{prob-revision}} & \textbf{RMSE} & \textbf{nRMSE} \\
\hline
0.03 & 0.145 & 20.6\% \\
\hline
0.06 & \textcolor{red}{0.120} & \textcolor{red}{17.1\%} \\
\hline
0.12 & 0.145 & 20.6\% \\
\hline
\end{tabular}
\end{table}

It is noteworthy that this conclusion is consistent with the optimal results observed in Case Study 1 (October 2024 dataset). As previously reported in Table~\ref{tab:prob_revision_comparison}, the best fitting performance on Dataset 1 was also achieved when \texttt{prob-revision} was set to 0.06 (RMSE = 0.114, nRMSE = 12.8\%). The consistency between the two sets of experimental results suggests that, under the current model architecture and behavioral mechanism, a strategy revision probability of 0.06 may represent the most ``reasonable'' or ``realistic'' social learning rhythm among individuals, enabling stable simulation of influence diffusion behaviors within the user population. Moreover, from the observed trends in the error metrics, it can be seen that setting \texttt{prob-revision} either too low (e.g., 0.03) or too high (e.g., 0.12) leads to a noticeable increase in model errors.

In further experiments, the errors between the model simulation results and the real-world topic evolution paths smoothed by the Moving Average (MA) method were calculated for Dataset 2. The results are summarized in Table~\ref{tab:smoothing_december_results}.

First, it was observed that regardless of the window size setting, prob-revision = 0.06 consistently achieved the lowest RMSE and nRMSE across all configurations. This finding further confirms that a revision probability of 0.06 exhibits strong fitting capability in the second dataset as well.

Second, as the smoothing window size increased from 4 to 7, a slight downward trend in both RMSE and nRMSE was observed across all \texttt{prob-revision} settings. This suggests that larger smoothing windows effectively reduce the interference of local fluctuations with overall trend matching, thereby slightly improving the model’s fitting performance.

Finally, \texttt{prob-revision} settings of 0.03 and 0.12 exhibited relatively similar error levels in this dataset, and both were significantly worse than 0.06. This result slightly differs from the observations in Case Study 1, where 0.03 also showed competitive performance. The discrepancy may be attributable to the faster information diffusion speed and more frequent individual responses in Dataset 2, making mid-frequency strategy revisions more aligned with the real-world evolution process.

\begin{table}[ht]
\centering
\caption{Simulation Errors under Different Moving Average Window Sizes (December Dataset)}
\label{tab:smoothing_december_results}
\begin{tabular}{|c|cc|cc|cc|}
\hline
\multirow{2}{*}{\textbf{\texttt{prob-revision}}} & \multicolumn{2}{c|}{\textbf{Window Size = 4}} & \multicolumn{2}{c|}{\textbf{Window Size = 5}} & \multicolumn{2}{c|}{\textbf{Window Size = 7}} \\
\cline{2-7}
 & RMSE & nRMSE & RMSE & nRMSE & RMSE & nRMSE \\
\hline
0.03 & 0.126 & 19.6\% & 0.124 & 19.4\% & 0.122 & 19.0\% \\
\hline
0.06 & 0.098 & 15.2\% & 0.096 & 15.0\% & \textcolor{red}{0.092} & \textcolor{red}{14.3\%} \\
\hline
0.12 & 0.126 & 19.6\% & 0.124 & 19.4\% & 0.121 & 18.9\% \\
\hline
\end{tabular}
\end{table}

Overall, across different smoothing scales, \texttt{prob-revision} = 0.06 consistently emerged as the optimal choice, and the application of smoothing (particularly with a window size of 7) further enhanced the model’s fitting performance relative to the real-world trends. These findings further reinforce the consistency and stability of the modeling framework when applied across different time periods and datasets, thereby validating the generalizability and robustness of the evolutionary modeling approach.

\subsubsection{Validation of Payoff Matrix Methods}

Using a fixed \texttt{prob-revision} value of 0.06, this section explores the performance of the simulations under different matrix construction methods. The results are then compared with those obtained from the first dataset (October 2024), providing a basis for cross-dataset analysis and evaluation.

\begin{table}[ht]
\centering
\caption{Simulation Errors under Different \(\alpha\) Values (December Dataset)}
\label{tab:alpha_december_results}
\begin{tabular}{|c|c|c|}
\hline
\textbf{\(\alpha\)} & \textbf{RMSE} & \textbf{nRMSE} \\
\hline
Baseline & \textcolor{orange}{0.120} & \textcolor{orange}{17.1\%} \\
\hline
0.1 & 0.120 & 17.1\% \\
\hline
0.2 & 0.136 & 19.5\% \\
\hline
0.3 & 0.143 & 20.4\% \\
\hline
0.4 & 0.128 & 18.3\% \\
\hline
0.5 & \textcolor{red}{0.108} & \textcolor{red}{15.4\%} \\
\hline
0.6 & \textcolor{blue}{0.152} & \textcolor{blue}{21.7\%} \\
\hline
0.7 & 0.138 & 19.7\% \\
\hline
0.8 & 0.138 & 19.8\% \\
\hline
0.9 & 0.143 & 20.4\% \\
\hline
\end{tabular}
\end{table}

The experimental results under different weighting coefficients \( \alpha \) are summarized in Table~\ref{tab:alpha_december_results}. Comparing the results with those from Dataset 1, several notable differences and commonalities were observed:

\begin{enumerate}
    \item \textbf{Similar Optimal \(\alpha\) Values:} In Dataset 1, the best simulation performance was achieved at \( \alpha = 0.3 \), whereas in Dataset 2, the optimal point shifted to \( \alpha = 0.5 \). This difference suggests that the underlying behavioral game mechanisms among users may vary structurally across different time periods or topical environments, highlighting the need for flexible parameter tuning to adapt the model to specific scenarios.

    \item \textbf{Same Error Trends:} In both datasets, model errors exhibited a noticeable increase when \( \alpha = 0.6 \), indicating that when agents' choice behaviors rely excessively on their own payoff, the resulting evolutionary trajectories are prone to deviating from real-world strategy dynamics.

    \item \textbf{Baseline Method Remains Valid:} Although tuning \( \alpha \) brought about performance improvements, the additive payoff structure (Baseline) consistently showed moderate performance across both datasets, confirming its value as a reasonable reference benchmark.
\end{enumerate}

In addition, further experiments were conducted to explore the fitting errors between the model simulation outputs and the real-world topic evolution paths after applying a Moving Average smoothing (window size = 7). Specifically, \texttt{prob-revision} was fixed at 0.06, and simulations were performed under different \( \alpha \) parameter settings. The corresponding RMSE and nRMSE results are summarized in Table~\ref{tab:alpha_smoothing_december_results}.

\begin{table}[ht]
\centering
\caption{Simulation Errors under Different \(\alpha\) Values (Smoothed Data, December Dataset)}
\label{tab:alpha_smoothing_december_results}
\begin{tabular}{|c|c|c|}
\hline
\textbf{\(\alpha\)} & \textbf{RMSE} & \textbf{nRMSE} \\
\hline
Baseline & \textcolor{orange}{0.092} & \textcolor{orange}{14.3\%} \\
\hline
0.1 & 0.093 & 14.5\% \\
\hline
0.2 & 0.113 & 17.6\% \\
\hline
0.3 & 0.119 & 18.6\% \\
\hline
0.4 & 0.103 & 16.0\% \\
\hline
0.5 & \textcolor{red}{0.078} & \textcolor{red}{12.1\%} \\
\hline
0.6 & \textcolor{blue}{0.130} & \textcolor{blue}{20.2\%} \\
\hline
0.7 & 0.114 & 17.7\% \\
\hline
0.8 & 0.114 & 17.7\% \\
\hline
0.9 & 0.119 & 18.5\% \\
\hline
\end{tabular}
\end{table}

From the results, it can be observed that the trend is consistent with the experiments conducted using unsmoothed raw data: the best fitting performance still appeared around \( \alpha = 0.5 \). Furthermore, the overall simulation errors (both RMSE and nRMSE) decreased after applying smoothing, leading to an improved model fitting to the real-world evolutionary trajectories.

\subsection{Impact of Interaction Dimension Weighting on Model Performance}

In the previous modeling processes, the influence score (payoff) of each information piece was calculated by assigning equal weights to the numbers of likes, comments, and reposts, namely through simple averaging. However, in real-world social media environments, different types of user interactions may contribute unequally to the formation of overall propagation influence. This section explores how different weighting settings across interaction dimensions affect the model's fitting performance, based on Dataset 1 (the October 2024 weather variation dataset).

First, the interaction dimension weights were set as \( \beta = 1 \), \( \gamma = 1.5 \), and \( \delta = 2.5 \), emphasizing the dominant role of reposts in the propagation process. Based on this configuration, the average payoff for each topic were recalculated. As shown in Table~\ref{tab:payoff_weighting_comparison}, \textit{Payoff1} refers to the original equal-weight setting (\( \beta = \gamma = \delta = 1 \)), while \textit{Payoff2} represents the payoff values after enhancing the weight of reposts (\( \delta = 2.5 \)) and moderately increasing the weight of comments (\( \gamma = 1.5 \)).

\begin{table}[ht]
\centering
\caption{Comparison of Topic Payoffs under Different Interaction Weightings}
\label{tab:payoff_weighting_comparison}
\begin{tabular}{|c|c|c|c|}
\hline
\textbf{Topic} & \textbf{Keyword} & \textbf{Payoff1} & \textbf{Payoff2} \\
\hline
0 & Temperature changes & 2.085529 & 3.254209 \\
\hline
1 & Air pollution & 1.518253 & 2.425714 \\
\hline
2 & Extreme weather & 2.110240 & 3.327851 \\
\hline
\end{tabular}
\end{table}

From the comparison results, it can be observed that the payoff values for all topics significantly increased after applying the weighted setting. Meanwhile, Topic 1, although its payoff value was also improved, consistently remained the least influential topic under both weighting schemes. This indicates that in this dataset, the discussion intensity surrounding air pollution was relatively weaker, and even when the weights of comments and reposts were increased, the overall attention advantage for this topic could not be significantly enhanced.

Under the above weighting configuration, further experiments were conducted to explore the optimal value for the \texttt{prob-revision} parameter. It was found that when \texttt{prob-revision} was set to 0.03, the simulation results achieved the lowest fitting error compared to the real data, indicating that a lower strategy revision probability better reflected the real evolutionary trends under the current weighting scheme. Subsequently, with \texttt{prob-revision} fixed at 0.03, systematic tests were performed across different \( \alpha \) parameter values to assess their impact on simulation performance.

\begin{figure}
    \centering
    \includegraphics[width=.9\linewidth]{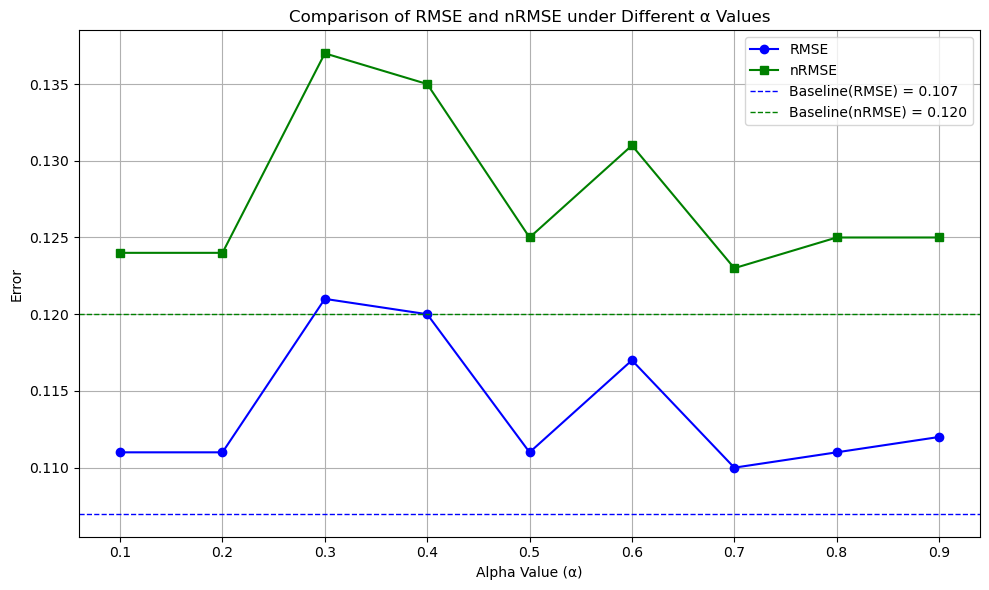}
    \caption{Visualisation of simulation Errors under Different \(\alpha\) Values (Original Data, \( \beta = 1 \), \( \gamma = 1.5 \), and \( \delta = 2.5 \))}
    \label{fig:1_15_25}
\end{figure}

As shown in Figure~\ref{fig:1_15_25}, although the overall fitting errors slightly increased and the simulation accuracy was marginally lower compared to the initial configuration with equal weights (\( \beta = \gamma = \delta = 1 \)), the differences in fitting errors across various \( \alpha \) values were within an acceptable range, suggesting that even after introducing interaction dimension weighting, the model was able to maintain reasonably good simulation performance.

\begin{figure}
    \centering
    \includegraphics[width=.9\linewidth]{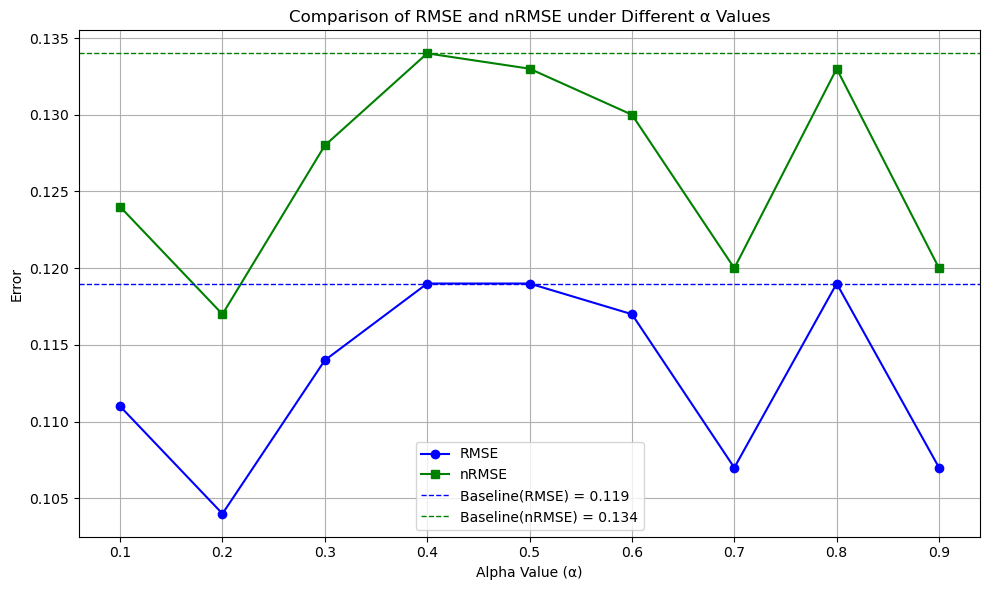}
    \caption{Visualisation of simulation Errors under Different \(\alpha\) Values (Original Data, \( \beta = 1 \), \( \gamma = 2.5 \), and \( \delta = 1.5 \))}
    \label{fig:1_25_15}
\end{figure}

Building on the previous exploration, this study further examined another set of interaction dimension weightings, specifically \( \beta = 1 \), \( \gamma = 2.5 \), and \( \delta = 1.5 \), to evaluate the impact of different assumptions regarding interaction importance on the model’s fitting performance. Preliminary experiments were first conducted to determine the optimal strategy revision probability under this weighting configuration. The results indicated that \texttt{prob-revision} = 0.03 again yielded the lowest model fitting error; thus, \texttt{prob-revision} was fixed at 0.03 for subsequent explorations of the \( \alpha \) parameter space. Under this weighting configuration, the simulation errors corresponding to different \( \alpha \) values are illustrated in Figure~\ref{fig:1_25_15}.

Based on the exploration results under the two sets of interaction weighting configurations, several findings can be summarized: Under the first parameter configuration (\( \beta = 1 \), \( \gamma = 1.5 \), \( \delta = 2.5 \)), although various \( \alpha \) values were tested, the overall model fitting performance did not exhibit significant improvements. The variations in RMSE and nRMSE were relatively small, and the best results were close to the baseline level. This suggests that under this setting, the constructed payoff structure already closely matched the diffusion characteristics of the subtopics, and further fine-tuning of the \( \alpha \) parameter provided limited additional benefits to the overall fitting accuracy. In contrast, under the second parameter configuration (\( \beta = 1 \), \( \gamma = 2.5 \), \( \delta = 1.5 \)), adjusting the \( \alpha \) parameter generally led to noticeable reductions in model errors. Multiple \( \alpha \) values achieved better fitting performance compared to the baseline setting. In particular, at \( \alpha = 0.2 \), \( \alpha = 0.7 \), and \( \alpha = 0.9 \), the model exhibited significant improvements in both RMSE and nRMSE metrics.

This phenomenon suggests that when the payoff weighting scheme changes, the relative weighting between self-perception and other-perception in individual payoff evaluation - controlled by the \( \alpha \) parameter
- exerts a nonlinear influence on the evolutionary trajectory. A preliminary hypothesis is that there may exist an upper threshold and a lower threshold within the model, determining whether adjusting \( \alpha \) can effectively improve system performance. Future research will further systematically explore the relationship between the payoff weighting configurations and \( \alpha \) sensitivity, aiming to investigate the existence of threshold regions and their potential impact on the dynamics of strategy evolution.

\section{Conclusion} \label{Conclusion}
This study has addressed a central challenge in computational social science: to move beyond isolated paradigms and build an integrative model that captures the dynamic interplay between the semantic content of online discourse and the behavioral mechanisms of its propagation. In response to the identified gap between content-centric and interaction-centric approaches, as well as the critique of idealized deliberative models, we have proposed and validated a novel computational framework that fuses advanced topic modeling, evolutionary game theory, and agent-based simulation.

Our work makes several key contributions. Theoretically, we have bridged the long-standing divide between semantic analysis and behavioral modeling by formally conceptualizing topics—extracted as coherent narratives from user-generated text—as strategic options in a social game. This allows us to model how the meaning of information directly influences its competitive fitness within the economy of attention, aligning with the concept of "cumulative deliberation" where public opinion forms through the accumulation of individual expressions. Methodologically, we have established a reproducible pipeline from raw textual data to a psychologically plausible simulation. By leveraging BERTopic \cite{grootendorst2022} for high-quality topic extraction and embedding evolutionary dynamics within a NetLogo \cite{tisue2004netlogo} ABM, we operationalize principles of bounded rationality and social learning \cite{bandura1971, rendell2010, simon1990}. Agents in our model make decisions through local imitation and payoff comparison, perturbed by behavioral noise, offering a more authentic representation of user behavior than models assuming perfect rationality. Empirically, our rigorous experiments on two real-world Weibo datasets concerning climate change have demonstrated the model's robustness and generalizability. We quantitatively showed how parameters like strategy revision probability (prob-revision) and payoff sensitivity ($\alpha$) control the system's evolution, providing nuanced insights into how micro-level imitation and cognitive constraints lead to emergent macro-level diffusion trends.

Despite its contributions, this study has several limitations that point to fruitful directions for future research. First, our model assumes a homogeneously mixed population, neglecting the complex, often scale-free, network structures that characterize real-world social platforms \cite{lerman2010, sasahara2021}. Integrating realistic network topologies is a crucial next step. Second, the agent's cognitive model is currently limited to a single dimension of topic preference. Future work should enrich the agent state space to include multi-dimensional factors such as emotional states \cite{fan2018}, pre-existing beliefs, social identity, and the effects of information overload \cite{kahneman2003}. Third, we have modeled topics as independent competitors. However, real-world discourses often exist in hierarchical, dependent, or symbiotic relationships. Exploring co-evolutionary topic dynamics and their interaction with platform algorithms represents a significant avenue for further study. Finally, the rapid advancement of LLMs promises even more powerful tools for semantic analysis. Future iterations of this framework could integrate generative models for dynamic topic summarization and real-time analysis, enhancing both the accuracy and timeliness of the simulations.

In summary, this research provides a foundational step towards a more holistic and interpretable modeling of complex social processes. By successfully integrating semantic analysis with a psychologically-grounded behavioral simulation, our framework not only offers a new pathway for social science inquiry powered by large language models but also delivers a practical tool for applications in misinformation resilience, opinion dynamics forecasting, and adaptive content curation. It underscores the critical importance of bridging computational methods with deeper theories of human cognition and social interaction to understand the ever-evolving landscape of online collective behavior.

\section{Acknowledgement} \label{Acknowledgement}
The research was supported by ...

\bibliographystyle{cas-model2-names}
\bibliography{cas-refs}

\end{document}